\documentstyle[12pt]{amsart}
\def	\t	{{\frak t}}
\def	\h	{{\frak h}}
\def	\R	{{\Bbb R}}
\def	\DH	{Duistermaat-Heckman\ }

\begin{document}
\title [Non-log-concave Duistermaat-Heckman measure]
{Example of a non-log-concave Duistermaat-Heckman measure}
\author{Yael Karshon}
\thanks{The author is partially supported by NSF grant DMS-9404404}
\thanks{dg-ga/9606007}
\address{Institute of Mathematics, The Hebrew University of Jerusalem,
Giv'at Ram, Jerusalem 91904, Israel}
\email{karshon@@math.huji.ac.il}
\begin{abstract}
We construct a compact symplectic manifold
with a Hamiltonian circle action for which the Duistermaat-Heckman
function is not log-concave.
\end{abstract}
\maketitle

\section{Introduction}
Let $T$ be a torus and $\t$ its Lie algebra.
Let $(M,\omega)$ be a symplectic manifold with an action of $T$
and with a moment map $$ \Phi : M \to \t^*.$$
Recall, this means that for every $\xi \in \t$, if $\xi_M$ is the 
corresponding vector field on $M$,  
$\iota(\xi_M) \omega = - d < \Phi, \xi >$.

{\bf Liouville measure} on $M$ associates to an open set $U$
the measure  $\int_U \omega^n$ where $n$ is half the dimension of
the manifold and where we integrate with respect to the symplectic
orientation.
The {\bf Duistermaat-Heckman measure} \cite{DH} on $\t^*$ 
is the push-forward of Liouville measure via the moment map $\Phi$.
If $T$ acts effectively, the Duistermaat-Heckman measure 
is absolutely continuous
with respect to Lebesgue measure, and the density function on $\t^*$ 
is called the {\bf Duistermaat-Heckman function}.

If $M$ is compact, the image of $\Phi$ is a convex polytope
\cite{GS,A}. If, in addition, the dimension of $T$ is half the
dimension of $M$ and $T$ acts effectively, the \DH function is 
equal to one on the convex polytope $\Phi(M)$ and zero outside 
\cite{Delzant}.  This function is log-concave, i.e., its logarithm
is concave.
Moreover, if we restrict this action to a subgroup $H$ of $T$, 
the moment map
for $H$ is the composition of the moment map for $T$ with the natural
linear projection  $\pi : \t^*  \to \h^*$. 
The \DH  function for $H$ is the function 
$x \mapsto {\operatorname{vol}} (\pi^{-1}(x) \cap \Phi(M))$
which associates to every point $x$ in $\h^*$
the volume of the corresponding ``slice" of the convex polytope $\Phi(M)$.
This function is again log-concave \cite[Theorem 6]{P}.

It was conjectured \cite{ginzburg,knudsen}
that for any Hamiltonian torus action on a compact manifold,
the \DH function is log-concave. This was proved for circle
actions on four manifolds in \cite[Remark 2.19]{karshon},
for coadjoint orbits in classical groups in \cite{okounkov},
and for arbitrary K\"ahler manifolds in \cite{graham}.
In this note we construct a counterexample to the conjecture;
we construct a Hamiltonian circle action on a compact symplectic manifold
for which \DH function is not log-concave. 
This construction came from investigating an example
of Dusa McDuff of a 6-manifold with a circle valued moment map 
\cite{mcduff}.  I use her notation wherever possible.

Our conventions regarding factors of $2\pi$ etc.\ are irrelevant
and will not be made explicit.

\subsection*{Acknowledgement}
I wish to thank Dror B-N.\ for commenting on the manuscript
and Y. Peres for advising me on convex sets and log-concave functions.

\section{The construction}
Let $T^4$ be the four dimensional torus with periodic coordinates
$x_i$, $1 \leq i \leq 4$, and let
$\sigma_{ij} = dx_i \wedge dx_j$ and 
$\sigma_{1234} = dx_1 \wedge dx_2 \wedge dx_3 \wedge dx_4$.
Let $L$ be a complex Hermitian line bundle over $T^4$
with Chern class $[ -\sigma_{14} - \sigma_{32} ]$.
Let $\Theta$ be a connection one-form with curvature 
$-\sigma_{14} - \sigma_{32}$.
This means that $\Theta$ is defined on $L$ outside the zero section,
that the restriction of $\Theta$ to a fiber of $L$ is $d\theta$
in polar coordinates on the fiber, and that 
$d\Theta$ is the pullback of $ - \sigma_{14} - \sigma_{32}$
via the bundle map $L \to T^4$. 
Denote by the same letters $\sigma_{ij},\sigma_{1234}$ the pullbacks
of these forms to $L$.
Let the function $\Phi : L \to \R$ be the norm squared, with respect
to the fiberwise Hermitian metric on $L$.  Consider the two-form
\begin{equation} \label{omega1}
  \omega = \sigma_{12} + \sigma_{34} 
  + (2-\Phi) \sigma_{14} + (3-\Phi) \sigma_{32} 
  + d\Phi \wedge \Theta
\end{equation}
on $L$ minus its zero section.
It is easy to check that $\omega$ is closed and that its top power is
$$\omega^3  =   6 (1 + (2-\Phi)(3-\Phi)) 
  \sigma_{1234} \wedge d\Phi \wedge \Theta .
$$
Since $\sigma_{1234} \wedge d\Phi \wedge \Theta \neq 0$
and since the function $(1 + (2-s)(3-s))$ is always positive,
$\omega$ is symplectic. 

The circle group acts on $L$ by fiberwise rotation. Let 
$\xi$ be the generating vector field.
From \eqref{omega1} it is clear that
$\iota(\xi) \omega = -d\Phi$, so $\Phi$ is a moment map for the 
circle action.
The \DH function is a constant positive multiple of the function
\begin{equation} \label{function}
\rho(s) = 1 + (2-s)(3-s).
\end{equation}
This function decreases for $0 < s < 2.5$ and increases for 
$2.5 < s < \infty$, so it is not log-concave.

To make a compact example out of our noncompact one, 
we perform ``Lerman cutting" \cite{lerman}: 
choose any two numbers, $0 < A < 2.5$
and $2.5 < B < \infty$.
``Lerman cutting" produces a compact symplectic manifold $(M,\omega)$
with a circle action and a moment map 
$\Phi : M \to [A,B]$ such that the preimages in $M$ and in $L$ 
of the open interval $(A,B)$ are equivariantly symplectomorphic.
Consequently, the \DH functions are the same: 
for the compact manifold $M$ we get the function \eqref{function}
restricted to the interval $A \leq s \leq B$, and this function
is not log-concave.

\end{document}